\documentclass[a4paper,10pt]{scrartcl}
\usepackage[english]{babel}
\usepackage[utf8]{inputenc}
\usepackage[T1]{fontenc}
\usepackage{lmodern} 
\usepackage{enumerate}
\usepackage[a4paper]{geometry}
\usepackage{amsfonts,amsmath,amssymb,amsthm,amsopn}
\usepackage{mathtools}
\usepackage{stmaryrd}
\usepackage{nicefrac}
\usepackage{exscale}

\usepackage[numbers,square]{natbib}
\usepackage{url} 
\usepackage[colorlinks=true,pdfpagelabels,unicode]{hyperref}
\hypersetup{urlcolor=red, citecolor=blue} 

\allowdisplaybreaks 

\newtheoremstyle{theo}
	{3pt} 
	{3pt} 
	{\itshape} 
	{} 
		{\bfseries} 
	{\\} 
	{ } 
	{\thmname{#1}\thmnumber{ #2.}\thmnote{ - #3}} 
\theoremstyle{theo}

\newtheorem{defi}{Definition}[section]
\newtheorem{lem}[defi]{Lemma}
\newtheorem{theo}[defi]{Theorem}
\newtheorem{cor}[defi]{Corollary}

\newtheorem{prop}[defi]{Proposition}

\newenvironment{bew}{\begin{proof}[\bfseries Proof:]}{\end{proof}}

\DeclareMathOperator{\bomega}{\overline{\Omega}}
\DeclareMathOperator{\romega}{\partial\Omega}

\DeclareMathOperator{\intd}{d\!}
\newcommand{\intdt}{\intd t}
\DeclareMathOperator{\dive}{\nabla\cdot}
\DeclareMathOperator{\divdot}{\cdot\!}

\newcommand{\Tm}{T_{max}}
\newcommand{\GNI}{Gagliardo-Nirenberg inequality}
\newcommand{\etd}{e^{t\Delta}}
\newcommand{\etA}{e^{-tA}}

\newcommand{\etsd}{e^{(t-s)\Delta}}

\newcommand{\intot}{\int_0^{t}} 

\newcommand{\intomega}{\int_{\Omega}\!} 

\newcommand{\intromega}{\int_{\romega}\!} 
\newcommand{\DA}[1][\delta]{D\!\left(A^{#1}\right)} 
\newcommand{\Lo}[1][1]{L^{#1}(\Omega)} 
\newcommand{\W}[1][1,2]{W^{#1}(\Omega)}
\newcommand{\Co}[1][0]{C^{#1}(\Omega)}
\newcommand{\LSp}[2]{L^{#1}\!\left(#2\right)} 

\newcommand{\CSp}[2]{C^{#1}\!\left(#2\right)}
\newcommand{\Po}{\mathcal{P}}
\newcommand{\nfrac}[2]{{\nicefrac{#1}{#2}}}

\newcounter{gleichung}
\newcommand{\owncount}{\refstepcounter{gleichung}}

\author{Tobias Black\thanks{Institut f\"ur Mathematik, Universit\"at Paderborn, Warburger Str. 100, 33098 Paderborn, Germany; email: \mbox{tblack@math.upb.de}}}
\title{Sublinear signal production in a two-dimensional Keller-Segel-Stokes system}

\begin{document}
\maketitle
\begin{abstract}
\noindent
{\textbf{Abstract:} 
We study the chemotaxis-fluid system 
\begin{align*}
\left\{\begin{array}{r@{\,}l@{\quad}l@{\,}c}
n_{t}\ &=\Delta n-\nabla\!\cdot(n\nabla c)-u\cdot\!\nabla n,\ &x\in\Omega,& t>0,\\
c_{t}\ &=\Delta c-c+f(n)-u\cdot\!\nabla c,\ &x\in\Omega,& t>0,\\
u_{t}\ &=\Delta u+\nabla P+n\cdot\!\nabla\phi,\ &x\in\Omega,& t>0,\\
\nabla\cdot u\ &=0,\ &x\in\Omega,& t>0,
\end{array}\right.
\end{align*}
where $\Omega\subset\mathbb{R}^2$ is a bounded and convex domain with smooth boundary, $\phi\in W^{1,\infty}\left(\Omega\right)$ and $f\in C^1([0,\infty))$ satisfies $0\leq f(s)\leq K_0 s^\alpha$ for all $s\in[0,\infty)$, with $K_0>0$ and $\alpha\in(0,1]$. This system models the chemotactic movement of actively communicating cells in slow moving liquid.

We will show that in the two-dimensional setting for any $\alpha\in(0,1)$ the classical solution to this Keller-Segel-Stokes-system is global and remains bounded for all times.
}
\end{abstract}

\section{Introduction}\label{sec1:intro}
\textbf{Keller-Segel models.}\quad Chemotaxis is the biological phenomenon of oriented movement of cells under influence of a chemical signal substance. This process is known to play a large role in various biological applications (\cite{HP09}). One of the first mathematical models concerning chemotaxis was introduced by Keller and Segel to describe the aggregation of bacteria (see \cite{KS70} and \cite{KS71}).
A simple realization of a standard Keller-Segel system, which models the assumption that the cells are not only attracted by higher concentration of the signal chemical but also produce the chemical themselves, can be expressed by
\begin{align}\label{KS}\tag{$KS$}
\left\{\begin{array}{r@{\,}l@{\quad}l@{\,}c}
n_{t}\ &=\Delta n-\nabla\!\cdot(n\nabla c),\ &x\in\Omega,& t>0,\\
c_{t}\ &=\Delta c-c+n,\ &x\in\Omega,& t>0,
\end{array}\right.
\end{align}
in a bounded domain $\Omega\subset\mathbb{R}^N$ with $N\geq1$. Herein, $n=n(x,t)$ denotes the unknown density of the involved cells and $c=c(x,t)$ the unknown concentration of the attracting chemical substance.

The Keller-Segel system alone has been studied intensively in the last decades and a wide array of interesting properties, such as  finite time blow-up and spatial pattern formation, have been discovered (see also the surveys \cite{BBWT15},\cite{HP09},\cite{Ho03}). For instance, the Keller-Segel system obtained from \eqref{KS} with homogeneous Neumann boundary conditions where $\Omega\subset\mathbb{R}^N$ is a ball, emits blow-up solutions for $N\geq2$, if the total initial mass of cells lies above a critical value (\cite{mizoguchi_winkler_13},\cite{win10jde}), while all solutions remain bounded when either $N=1$, or $N=2$ an the initial total mass of cells is below the critical value (\cite{OY01},\cite{NSY97}).

Through its application to various biological contexts, many variants of the Keller-Segel model have been proposed over the  years. In particular, adaptions of \eqref{KS} in the form of
\begin{align}\label{kssens}
\owncount
n_{t} =\Delta n-\nabla\!\cdot(nS(x,n,c)\cdot\nabla c),\quad x\in\Omega, t>0,
\end{align}
with given chemotactic sensitivity function $S$, which can either be a scalar function, or more general a tensor valued function (see e.g. \cite{XO09-MSmodels}), for the first equation or
\begin{align}\label{ksox}
\owncount
c_{t} =\Delta c-ng(c),\quad x\in\Omega, t>0,
\end{align}
with given function $g$ for the second equation, have been studied. Both of these adjustments are known to have an influence on the boundedness of solutions to their respective systems. For instance, if we replace the first equation of \eqref{KS} with \eqref{kssens} for a scalar function $S$ satisfying $S(r)\leq C(1+r)^{-\gamma}$ for all $r\geq1$ and some $\gamma>1-\frac{2}{N}$, then all solutions to the corresponding Neumann problem are global and uniformly bounded. On the other hand if $N\geq2$, $\Omega\subset\mathbb{R}^N$ is a ball and $S(r)>cr^{-\gamma}$ for some $\gamma<1-\frac{2}{N}$ then the solution may blow up (\cite{HoWin05_bvblowchemo}). 

Considering the adaption of \eqref{KS} with \eqref{ksox} as second equation, which basically corresponds to the system assumption that the cells consume some of the chemical instead of producing it, it was shown in \cite{TaoWin12_evsmooth} that for $N=2$ the corresponding Neumann problem possesses bounded classical solution for suitable regular initial data not depending on a smallness condition. For $N=3$ it was proved, that there exist global weak solutions which eventually become smooth and bounded after some waiting time. 

A combination of both adjustments, where $S$ is matrix-valued with non-trivial nondiagonal parts, was studied in \cite{win15_chemorot}. There it was shown that under fairly general assumptions on $g$ and $S$ at least one generalized solution exists which is global. This result does neither contain a restriction on the spatial dimension nor on the size of the initial data.

One last adaption of \eqref{KS} we would like to mention has only recently been studied thoroughly and concerns the system
\begin{align}\label{KSa}\tag{$KS^\alpha$}
\left\{\begin{array}{r@{\,}l@{\quad}l@{\,}c}
n_{t}\ &=\Delta n-\nabla\!\cdot(n\nabla c),\ &x\in\Omega,& t>0,\\
c_{t}\ &=\Delta c-c+f(n),\ &x\in\Omega,& t>0,
\end{array}\right.
\end{align}
with $f\in C^1\left([0,\infty)\right)$ satisfying $0\leq f(n)\leq K n^\alpha$ for any $n\geq0$ with $K>0$ and $\alpha>0$. In this setting it is known, that the system \eqref{KSa} does not emit any blow-up solution if $\alpha<\frac{2}{N}$ (\cite{liudongmei15_boundchemo}) but it remains an open question whether this exponent is indeed critical.

Similar forms of $f(n)$ have been treated before either in the linear case $f(n)=n$ (\cite{Mimura1996499}) or (sub-)linear cases with an additional logistic growth term introduced to the first equation (eg. \cite{Os02-chemologatract},\cite{Win10-chemolog},\cite{NaOs13}).

\textbf{Chemotaxis-fluid systems.}\quad Nonetheless, one assumption is shared by all of these adapted Keller-Segel models. That is, only the cell density $n$ and the chemical concentration $c$ are unknown and all other system parameters are fixed. In particular, the models assume that there is no interaction between the cells and their surroundings. However, experimental observations indicate that chemotactic motion inside a liquid can be substantially influenced by the mutual interaction between cells and fluid. For instance, in \cite{tuval2005bacterial} the dynamical generation of patterns and emergence of turbulence in population of aerobic bacteria suspended in sessile drops of water is reported, whereas examples involving instationary fluids are important in the context of broadcast spawning phenomena related to successful coral fertilization (\cite{coll1994chemical},\cite{miller1985demonstration}).

A model considering the chemotaxis-fluid interaction, building on experimental observations of Bacillus subtilis was given in \cite{tuval2005bacterial}. In the system in question, the fluid velocity $u=u(x,t)$ and the associated pressure $P=P(x,t)$ are introduced as additional unknown quantities utilizing the incompressible Navier-Stokes equations. One of the first theoretical results concerning the solvability in this context were shown in \cite{lorz10}, where the local existence of weak solutions for $N\in\{2,3\}$ was shown. This setting, however, involved signal consumption in the form of per-capita oxygen consumption of the bacteria, which corresponds to an equation of the form \eqref{ksox}. Since we want to focus on the case of signal production by the cells as realized in \eqref{KS}, a more suitable system in this context is the Keller-Segel-Navier-Stokes system
\begin{align}\label{KSNS}\tag{$KSNS$}
\left\{\begin{array}{r@{\,}r@{\,}r@{\,}l@{\quad}l@{\,}c}
n_{t}\,&+\,&u\cdot\!\nabla n\ &=\Delta n-\nabla\!\cdot(n\nabla c),\ &x\in\Omega,& t>0,\\
c_{t}\,&+\,&u\cdot\!\nabla c\ &=\Delta c-c+n,\ &x\in\Omega,& t>0,\\
u_{t}\,&+\,&u\cdot\!\nabla u\ &=\Delta u-\nabla P+n\cdot\!\nabla\phi,\ &x\in\Omega,& t>0,\\
&&\dive u\ &=0,\ &x\in\Omega,& t>0,
\end{array}\right.
\end{align}
where the fluid is supposed to be driven by forces induced by the fixed gravitational potential $\phi$ and transports both the cells and the chemical.

The mathematical analysis of \eqref{KSNS} regarding global and bounded solutions is far from trivial, as on the one hand its Navier-Stokes subsystem lacks complete existence theory (\cite{Wie99-NS}) and on the other hand the previously mentioned properties for Keller-Segel system can still emerge. In order to weaken the analytical effort necessary, a commonly made simplification is to assume that the fluid flow is comparatively slow and thus the fluid velocity evolution may be described by the Stokes equation $(\kappa=0)$ rather than the full Navier-Stokes system.

Of course, all alterations to \eqref{KS} described above can be included as adjustments to the systems in this Keller-Segel(-Navier)-Stokes setting as well. Their influence on global and bounded solutions are one focal point of recent studies. For instance, an adjustment making use of both sensitivity and chemical consumption has been applied to Keller-Segel-Stokes systems in \cite{win15_globweak3d}, where for scalar valued sensitivity functions $S$ the existence of global weak solutions for bounded three-dimensional domains has been established. Building on this existence result, it was shown in \cite{win15_chemonavstokesfinal} that the generalized solution approaches a spatially homogeneous steady state under fairly weak assumptions imposed on the parameter functions $S$ and $g$. Under similar assumptions the existence of global weak solutions for suitable non-linear diffusion types have been proven in \cite{francescolorz10} and the existence of bounded and global weak solutions even allowing matrix-valued $S$ not requiring a decay assumption in \cite{win_ct_fluid_3d}.

A Keller-Segel-Stokes system corresponding to the adjustment made to \eqref{KS} by only making use of rotational sensitivity was studied in \cite{Wang20157578}, where it was shown that the Neumann problem for the Keller-Segel-Stokes system possesses a unique global classical solution which remains bounded for all times, if we assume $S$ to satisfy $|S(x,n,c)|\leq C_S(1+n)^{-\alpha}$ with $C_S>0$ for some $\alpha>0$.

Regarding the introduction of the additional logistic growth term $+rn-\mu n^2$ with $r\geq0$ and $\mu>0$ to the first equation, it was shown in \cite[Theorem 1.1]{tao_winkler15_zampfinal} for space dimension $N=3$, that every solution remains bounded if $\mu\geq23$ and thus any blow-up phenomena are excluded. Moreover, these solutions tend to zero (\cite[Theorem 1.2]{tao_winkler15_zampfinal}).

Some of these results have in part been transferred to the full chemotaxis Navier-Stokes system. These include global existence of classical solutions for $N=2$ with scalar valued sensitivity (\cite{win_fluid_final}), large time behavior and eventual smoothness of such solutions (\cite{win15_chemonavstokesfinal}) and even global existence of mild solution to double chemotaxis systems under the effect of incompressible viscious fluid (\cite{kozono15}). Boundedness results with matrix-valued sensitivity without decay requirements but for small initial data have been discussed in \cite{caolan16_smalldatasol3dnavstokes} and boundedness results under influence of a logistic growth term in \cite{tao_winkler_non2015}.

\textbf{Main results.}\quad  Most of these results stated above are concerned with the chemical consumption version of the chemotaxis model (\cite{Wang20157578} and \cite{tao_winkler_non2015} being the exceptions). To the best of our knowledge the Stokes variant of chemotaxis-fluid interaction has only been discussed outside of the chemical consumption case  either by introducing a logistic growth term as in \cite{tao_winkler_non2015} or taking a more general chemotactic sensitivity as in \cite{Wang20157578}. Motivated by this fact and the result of \cite{liudongmei15_boundchemo} for \eqref{KSa} mentioned above, we are now interested in whether the influence of a coupled slow moving fluid described by Stokes equation affects the possible choice for $\alpha\in(0,1)$, while still maintaining the exclusion of possible unbounded solutions. Henceforth, we will consider that the evolution of $(n,c,u,P)$ is governed by the Keller-Segel-Stokes System
\begin{align}\label{KSaS}\tag{$KS^{\alpha}S$}
\left\{\begin{array}{r@{\,}l@{\quad}l@{\,}c}
n_{t}\ &=\Delta n-\nabla\!\cdot(n\nabla c)-u\cdot\!\nabla n,\ &x\in\Omega,& t>0,\\
c_{t}\ &=\Delta c-c+f(n)-u\cdot\!\nabla c,\ &x\in\Omega,& t>0,\\
u_{t}\ &=\Delta u+\nabla P+n\cdot\!\nabla\phi,\ &x\in\Omega,& t>0,\\
\dive u\ &=0,\ &x\in\Omega,& t>0,
\end{array}\right.
\end{align}
where $\Omega\subset\mathbb{R}^2$ is a bounded and smooth domain and $f\in C^1([0,\infty))$ satisfies 
\begin{align}\label{fprop}
\owncount
0\leq f(s)\leq K_0 s^\alpha\quad\mbox{for all }s\in[0,\infty)
\end{align}
with some $\alpha\in(0,1]$ and $K_0>0$. We shall examine this system along with no-flux boundary conditions for $n$ and $c$ an a no-slip boundary condition for $u$,
\begin{align}\label{bcond}
\owncount
\frac{\partial n}{\partial\nu}=\frac{\partial c}{\partial\nu}=0\quad\mbox{and}\quad u=0\qquad \mbox{for }x\in\romega\mbox{ and }t>0,    
\end{align}
and initial conditions
\begin{align}\label{idcond}
\owncount
n(x,0)=n_0(x),\quad c(x,0)=c_0(c),\quad u(x,0)=u_0(x),\quad x\in\Omega.
\end{align}
For simplicity we will assume $\phi\in\W[1,\infty]$ and that for some $\theta>2$ and $\delta\in(\frac{1}{2},1)$ the initial data satisfy the regularity and positivity conditions
\begin{align}\label{idreg}
\owncount
\begin{cases}&n_0\in \CSp{0}{\bomega}\mbox{ with }n_0> 0\mbox{ in }\bomega,\\
&c_0\in\W[1,\theta]\mbox{ with }c_0> 0\mbox{ in }\bomega,\\
&u_0\in \DA,\end{cases}
\end{align}
where here and below $A^\delta$ denotes the fractional power of the Stokes operator $A:=-\Po\Delta$ in $\Lo[2]$ regarding homogeneous Dirichlet boundary conditions, with the Helmholtz projection $\Po$ from $\Lo[2]$ to the solenodial subspace $L^2_\sigma(\Omega):=\left\{\left.\varphi\in\Lo[2]\right\vert\dive\varphi=0\right\}$.

In this framework we can state our main result in the following way:
\begin{theo}\label{Thm:globEx}
Let $\theta>2$, $\delta\in(\frac{1}{2},1)$ and $\Omega\subset\mathbb{R}^2$ be a bounded and convex domain with smooth boundary. Assume $\phi\in\W[1,\infty]$ and that $n_0,c_0$ and $u_0$ comply with \eqref{idreg}. Then for any $\alpha\in(0,1)$, the PDE system \eqref{KSaS} coupled with boundary conditions \eqref{bcond} and initial conditions \eqref{idcond} possesses a solution $(n,c,u,P)$ satisfying
\begin{align*}
\begin{cases}
n\in\CSp{0}{\bomega\times[0,\infty)}\cap\CSp{2,1}{\bomega\times(0,\infty)},\\
c\in\CSp{0}{\bomega\times[0,\infty)}\cap\CSp{2,1}{\bomega\times(0,\infty)},\\
u\in\CSp{0}{\bomega\times[0,\infty)}\cap\CSp{2,1}{\bomega\times(0,\infty)},\\
P\in\CSp{1,0}{\bomega\times[0,\infty)},
\end{cases}
\end{align*}
which solves \eqref{KSaS} in the classical sense and remains bounded for all times. This solution is unique within the class of functions which for all $T\in(0,\infty)$ satisfy the regularity properties
\begin{align}\label{locExUClass}
\owncount
\begin{cases}
n\in\CSp{0}{[0,T);\Lo[2]}\cap\LSp{\infty}{(0,T);\CSp{0}{\bomega}}\cap\CSp{2,1}{\bomega\times(0,T)},\\
c\in\CSp{0}{[0,T);\Lo[2]}\cap\LSp{\infty}{(0,T);\W[1,\theta]}\cap\CSp{2,1}{\bomega\times(0,T)},\\
u\in\CSp{0}{[0,T);\Lo[2]}\cap\LSp{\infty}{(0,T);\DA}\cap\CSp{2,1}{\bomega\times(0,T)},\\
P\in \LSp{1}{(0,T);\W[1,2]},
\end{cases}
\end{align}
up to addition of functions $\hat{p}$ to $P$, such that $\hat{p}(\cdot,t)$ is constant for any $t\in(0,\infty)$.
\end{theo}
In view of Theorem \ref{Thm:globEx}, there is no evident difference regarding $\alpha$ between the coupled system \eqref{KSaS} and the chemotaxis system without fluid \eqref{KSa} for dimension $N=2$.

In Section \ref{sec1:lEx} we will briefly discuss local existence of classical solutions and basic a priori estimates. Section \ref{sec2:regu} is dedicated to the connection between the regularity of $n$ and the regularity of the spacial derivative of $u$, which plays a crucial part in obtaining additional information on the regularity of $c$. In Section \ref{sec3:N2} we will combine standard testing procedures with the results from the previous sections to prove the boundedness and globality of classical solutions to \eqref{KSaS}.

%
\setcounter{gleichung}{0} 
\section{Local existence of classical solutions}\label{sec1:lEx}
The following theorem concerning the local existence of classical solutions, as well as an extensibility criterion can be proven with exactly the same steps demonstrated in \cite[Lemma 2.1]{win_fluid_final} and \cite[Lemma 2.1]{tao_winkler_chemohapto11siam}.

\begin{lem}[Local existence of classical solutions]\label{Lem:locEx}
Let $\theta>2$, $\delta\in(\frac{1}{2},1)$ and $\Omega\subset\mathbb{R}^2$ be a bounded and convex domain with smooth boundary. Suppose $\phi\in\W[1,\infty]$ and that $n_0,c_0$ and $u_0$ satisfy \eqref{idreg}. Then there exist $\Tm\in(0,\infty]$ and functions $(n,c,u,P)$ satisfying
\begin{align*}
\begin{cases}
n\in\CSp{0}{\bomega\times[0,\Tm)}\cap\CSp{2,1}{\bomega\times(0,\Tm)},\\
c\in\CSp{0}{\bomega\times[0,\Tm)}\cap\CSp{2,1}{\bomega\times(0,\Tm)},\\
u\in\CSp{0}{\bomega\times[0,\Tm)}\cap\CSp{2,1}{\bomega\times(0,\Tm)},\\
P\in\CSp{1,0}{\bomega\times[0,\Tm)},
\end{cases}
\end{align*}
that solve \eqref{KSaS} with \eqref{bcond} and \eqref{idcond} in the classical sense in $\Omega\times(0,\Tm)$. Moreover, we have $n>0$ and $c>0$ in $\bomega\times[0,\Tm)$ and the alternative
\begin{align}\label{lExAlt}
\owncount
\mbox{either }\,\Tm=\infty\,\mbox{ or }\,\|n(\cdot,t)\|_{\Lo[\infty]}+\|c(\cdot,t)\|_{\W[1,\theta]}+\|A^{\delta}u(\cdot,t)\|_{\Lo[2]}\to\infty\ \mbox{as }t\nearrow\Tm.
\end{align}
This solution is unique among all functions satisfying \eqref{locExUClass} for all $T\in(0,\Tm)$, up to addition of functions $\hat{p}$, such that $\hat{p}(\cdot,t)$ is constant for any $t\in(0,T)$, to the pressure $P$.
\end{lem}

Local existence at hand, we can immediately prove two elementary properties, which will be the starting point for all of our regularity results to come.

\begin{lem}\label{Lem:masscons}
Under the assumptions of Lemma \ref{Lem:locEx}, the solution of \eqref{KSaS} satisfies
\begin{align}\label{masscons-n}
\owncount
\intomega n(x,t)\intd x=\intomega n_0=:m\quad\mbox{for all}\ t\in(0,\Tm)
\end{align}
and there exists a constant $C>0$ such that
\begin{align}\label{L1-bound-c}
\owncount
\intomega c(x,t)\intd x\leq C\quad\mbox{for all}\ t\in(0,\Tm).
\end{align}
\end{lem}

\begin{bew}
The first property follows immediately from simple integration of the first equation in \eqref{KSaS}. For \eqref{L1-bound-c} we integrate the second equation of \eqref{KSaS} and recall \eqref{fprop} to obtain
\begin{align*}
\frac{\intd}{\intdt}\intomega c+\intomega c\leq K_0\intomega n^\alpha\quad\mbox{for all }t\in(0,\Tm).
\end{align*}
Hence, making use of \eqref{masscons-n} and the fact $\alpha<1$, $y(t)=\intomega c(x,t)\intd x$ satisfies the ODI
\begin{align*}
y'(t)+y(t)\leq C_1\|n_0\|_{\Lo[1]}^{\alpha}=C_2\quad\mbox{for all }t\in(0,\Tm)
\end{align*}
for some $C_1>0$ and $C_2:=C_1m^\alpha>0$ in view of \eqref{idreg}. Upon integration we infer 
\begin{align*}
y(t)\leq y(0)e^{-t}+C_2\left(1-e^{-t}\right)\quad\mbox{for all }t\in(0,\Tm),
\end{align*}
which, due to the assumed regularity of $c_0$ in \eqref{idreg}, completes the proof.
\end{bew}

\setcounter{gleichung}{0} 
\section{Regularity of u implied by regularity of n}\label{sec2:regu}
Let us recall that $\Po$ denotes the Helmholtz projection from $\Lo[2]$ to the subspace $L^2_\sigma\left(\Omega\right)=\left\{\varphi\in\Lo[2]\,\vert\,\dive\varphi=0\right\}$ and $A:=-\Po\Delta$ denotes the Stokes operator in $\Lo[2]$ under homogeneous Dirichlet boundary conditions.

For now we limit our observations to a projected version of the Stokes subsystem $\frac{\intd}{\intd t}u+Au=\Po\left(n\nabla\phi\right)$ in \eqref{KSaS} without regard for the rest of the system. In contrast to the setting with the full Navier-Stokes equations we can make use of the absence of the convective term $(u\cdot\nabla)u$ in the Stokes equation to gain results concerning the regularity of the spatial derivative $Du$ based on the regularity of the term $\Po\left(n\nabla\phi\right)$, which in fact solely depends on the regularity of $n$, due to the assumed boundedness of $\nabla\phi$.

In \cite[Lemma 2.4]{Wang20157578} this correlation between the regularity of $u$ and $n$ is proven in space dimension $N=2$. The proof of \cite[Lemma 2.4]{Wang20157578} is based on an approach employed in \cite[Section 3.1]{win_ct_fluid_3d}, which makes use of general results for sectorial operators shown in \cite{fr69}, \cite{hen81} and \cite{gig81} and mainly relies on an embedding of the domains of fractional powers $D\left(A^\beta\right)$ into $\Lo[p]$, see \cite[Theorem 1.6.1]{hen81} or \cite[Theorem 3]{gig81}, for instance. Since we are only working in two-dimensional domains we will only state the result from \cite[Lemma 2.4]{Wang20157578} here and refer the reader to \cite[Corollary 3.4]{win_ct_fluid_3d} and \cite[Lemma 2.5]{Wang20157578} for the remaining details regarding the proof.

\begin{lem}\label{Lem:u_w1r_from_n_lp}
Let $p\in[1,\infty)$ and $r\in[1,\infty]$ be such that
\begin{align*}
\begin{cases}
r<\frac{2p}{2-p}\quad&\mbox{if }p\leq 2,\\
r\leq\infty\quad&\mbox{if }p>2.
\end{cases}
\end{align*}
Furthermore, let $T>0$ be such that $n:\Omega\times(0,T)\mapsto\mathbb{R}$ satisfies
\begin{align*}
\|n(\cdot,t)\|_{\Lo[p]}\leq \eta\quad\mbox{for all }t\in(0,T),
\end{align*}
with some $\eta>0$. Then for $u_0\in\DA$ with $\delta\in\left(\frac{1}{2},1\right)$ and $\phi\in\W[1,\infty]$ all solutions $u$ of the third and fourth equations in \eqref{KSaS} fulfill
\begin{align*}
\|Du(\cdot,t)\|_{\Lo[r]}\leq C\quad\mbox{for all }t\in(0,T),
\end{align*}
with a constant $C=C(p,r,\eta,u_0,\phi)>0$.
\end{lem}

Evidently, a supposedly known bound for $n$ at hand, we immediately obtain the desired boundedness of $u$ in view of Sobolev embeddings. Nevertheless, since we only have the time independent $L^1$--bound of $n$ from Lemma \ref{Lem:masscons} as a starting point, obtaining a bound for $n$ in $\Lo[p]$ with suitable large $p>1$ will require additional work.

\setcounter{gleichung}{0} 
\section{Global existence and boundedness in two-dimensional domains}\label{sec3:N2}
For this section we fix $\theta>2,\delta\in(\frac{1}{2},1)$ and initial data satisfying \eqref{idreg}. In particular, this ensures that all requirements of Lemma \ref{Lem:locEx} are met. Let $(n,c,u,P)$ denote the solution given by Lemma \ref{Lem:locEx} and $\Tm$ its maximal time of existence. Making use of the connection between the regularity of $u$ and $n$ discussed in the previous section we immediately obtain

\begin{prop}\label{Prop:u_bound_N2_all_p}
For all $r<2$ and all $q<\infty$ there exist constants $C_1>0$ and $C_2>0$ such that the solution to \eqref{KSaS} satisfies
\begin{align*}
\|Du(\cdot,t)\|_{\Lo[r]}\leq C_1\quad\mbox{for all }t\in(0,\Tm)
\end{align*}
and
\begin{align*}
\|u(\cdot,t)\|_{\Lo[q]}\leq C_2\quad\mbox{for all }t\in(0,\Tm).
\end{align*}
\end{prop}

\begin{bew}
Due to \eqref{masscons-n} and the regularity of $n_0$ we can find $C_3>0$ such that $\|n(\cdot,t)\|_{\Lo[1]}=\|n_0\|_{\Lo[1]}\leq C_3$ holds for all $t\in(0,\Tm)$. Thus, we may apply Lemma \ref{Lem:u_w1r_from_n_lp} with $p=1$ to obtain for any $r<2$ that $\|Du(\cdot,t)\|_{\Lo[r]}\leq C_2\ \mbox{for all }t\in(0,\Tm)$ with some $C_2>0$. The second claim then follows immediately from the Sobolev embedding theorem (\cite[Theorem 5.6.6]{evans}).
\end{bew}

\subsection{Obtaining a first information on the gradient of c}\label{ssec32:reg_c}

In order to derive the bounds necessary in our approach towards the boundedness result we require an estimate on the gradient of $c$ as a starting point. To obtain a first information in this matter, we apply standard testing procedures to derive an energy inequality involving integrals of $n\ln n$ and $|\nabla c|^2$. But first, let us briefly recall Young's inequality in order to fix notation.

\begin{lem}\label{young}
Let $a,b,\varepsilon>0$ and $1<p,q<\infty$ with $\frac{1}{p}+\frac{1}{q}=1$. Then
\begin{align*}
ab\leq\varepsilon a^p+C(\varepsilon,p,q) b^q,
\end{align*}
where $C(\varepsilon,p,q)=(\varepsilon p)^{-\frac{q}{p}}q^{-1}$.
\end{lem}

Before we derive an inequality for the time evolution of $\intomega n\ln$ we employ the \GNI\ to show one simple preparatory lemma on which we will rely multiple times later on.

\begin{lem}\label{Lem:ngnb}
Let $\Omega\subset\mathbb{R}^2$ be a bounded domain with smooth boundary. Let $r\geq1$ and $s\geq1$. Then for any $\eta>0$ there exists $C>0$ such that
\begin{align*}
\intomega|\varphi|^{rs}\leq C\left(\intomega|\nabla(|\varphi|^\nfrac{r}{2})|^2\right)^{\frac{(rs-1)}{r}}+C
\end{align*}
holds for all functions $\varphi\in\Lo[1]$ satisfying $\nabla(|\varphi|^{\nfrac{r}{2}})\in\Lo[2]$ and $\intomega|\varphi|\leq\eta$.
\end{lem}

\begin{bew}
By an application of the \GNI\ (see \cite[Lemma 2.3]{lankchapto15} for a version including integrability exponents less than $1$) we can pick $C_1>0$ such that
\begin{align*}
\intomega|\varphi|^{rs}=\||\varphi|^\nfrac{r}{2}\|_{\Lo[2s]}^{2s}\leq C_1\|\nabla(|\varphi|^\nfrac{r}{2})\|_{\Lo[2]}^{2sa}\||\varphi|^\nfrac{r}{2}\|_{\Lo[\nfrac{2}{r}]}^{2s(1-a)}+C_1\||\varphi|^\nfrac{r}{2}\|_{\Lo[\nfrac{2}{r}]}^{2s}
\end{align*}
holds for all $\varphi\in\Lo[1]$ with $\nabla(|\varphi|^{\nfrac{r}{2}})\in\Lo[2]$, with $a\in(0,1)$ provided by
\begin{align*}
a=\frac{\frac{r}{2}-\frac{1}{2s}}{\frac{r}{2}+\frac{1}{2}-\frac{1}{2}}=1-\frac{1}{rs}.
\end{align*}
Since $\intomega|\varphi|\leq\eta$ we have $\||\varphi|^\nfrac{r}{2}\|_{\Lo[\nfrac{2}{r}]}=\left(\intomega|\varphi|\right)^\nfrac{r}{2}\leq\eta^\nfrac{r}{2}$ and thus
\begin{align*}
\intomega|\varphi|^{rs}\leq C_2\left(\intomega|\nabla(|\varphi|^\nfrac{r}{2})|^2\right)^{\frac{(rs-1)}{r}}+C_2
\end{align*}
for all $\varphi\in\Lo[1]$ satisfying $\nabla(|\varphi|^{\nfrac{r}{2}})\in\Lo[2]$, where $C_2=C_1\max\{\eta,\eta^{rs}\}>0$.
\end{bew}

The particular form in which we will need this inequality most often is the following:

\begin{cor}\label{cor:ngnb}
There exists a constant $K_1>0$ such that the solution of \eqref{KSaS} fulfills
\begin{align*}
\intomega n^2\leq K_1\intomega|\nabla(n^\nfrac{1}{2})|^2+K_1
\end{align*}
for all $t\in(0,\Tm)$.
\end{cor}

Testing of the first equation in \eqref{KSaS} with $\ln n$ yields the following estimation.

\begin{lem}\label{Lem:n-energy}
There exists a constant $K_2>0$ such that the solution of \eqref{KSaS} fulfills
\begin{align}\label{n-energy}
\owncount
\frac{\intd}{\intdt}\intomega n\ln n +\intomega|\nabla( n^{\nfrac{1}{2}})|^2\leq K_2\intomega\left|\Delta c\right|^2+K_2\quad\mbox{for all }t\in(0,\Tm).
\end{align}
\end{lem}

\begin{bew}
Making use of \eqref{masscons-n} and $\dive u=0$ in $\Omega$, multiplication of the first equation in \eqref{KSaS} with $\ln n$ and integration by parts yield
\begin{align}\label{n-energy-proof-eq1}
\owncount
\frac{\intd}{\intdt}\intomega n\ln n+\intomega\frac{|\nabla n|^2}{n}=\intomega\nabla c\divdot\nabla n\quad\mbox{for all }t\in(0,\Tm).
\end{align}
To further estimate the right hand side, we first let $K_1>0$ be as in Corollary \ref{cor:ngnb}. Then, integrating the right hand side of \eqref{n-energy-proof-eq1} by parts once more and applying Young's inequality with $p=q=2$ and $\varepsilon=\frac{3}{K_1}$ (see Lemma \ref{young}) and Corollary \ref{cor:ngnb}, we obtain
\begin{align*}
\frac{\intd}{\intdt}\intomega n\ln n+4\intomega|\nabla(n^{\nfrac{1}{2}})|^2&\leq\frac{3}{K_1}\intomega n^2+C_1\intomega|\Delta c|^2\\
&\leq\frac{3}{K_1}\left(K_1\intomega|\nabla(n^{\nfrac{1}{2}})|^2+K_1\right)+C_1\intomega|\Delta c|^2
\end{align*}
for all $t\in(0,\Tm)$ and some $C_1>0$. Reordering the terms appropriately completes the proof with $K_2:=\max\{3,C_1\}$.
\end{bew}

The second of the separate inequalities treats the time evolution of $\intomega|\nabla c|^2$.

\begin{lem}\label{Lem:c-energy}
Given any $\xi>0$, there exists a constant $K_3>0$ such that
\begin{align}\label{c-energy}
\owncount
\frac{\xi}{2}\frac{\intd}{\intdt}\intomega |\nabla c|^2 +\frac{\xi}{4}\intomega|\Delta c|^2+ \xi\intomega|\nabla c|^2\leq \frac{1}{2}\intomega|\nabla(n^{\nfrac{1}{2}})|^2+K_3
\end{align}
holds for all $t\in(0,\Tm)$.
\end{lem}

\begin{bew}
Testing the second equation of \eqref{KSaS} with $-\xi\Delta c$ and integrating by parts we obtain 
\begin{align*}
\frac{\xi}{2}\frac{\intd}{\intdt}\intomega|\nabla c|^2 +\xi\intomega|\Delta c|^2+\xi\intomega|\nabla c|^2=-\xi\intomega f(n)\Delta c+\xi\intomega\Delta c\nabla c\cdot u
\end{align*}
for all $t\in(0,\Tm)$. An application of Young's inequality to both integrals on the right side therefore implies that
\begin{align}\label{c-energy-eq0}
\owncount
\frac{\xi}{2}\frac{\intd}{\intdt}\intomega|\nabla c|^2+\xi\intomega|\Delta c|^2+\xi\intomega|\nabla c|^2 \leq \xi\intomega f(n)^{2}+\frac{\xi}{2}\intomega|\Delta c|^2+\xi\intomega|\nabla c|^2 |u|^2
\end{align}
holds for all $t\in(0,\Tm)$. We fix $q>2$ and make use of the Hölder inequality to see that
\begin{align}\label{c-energy-eq1}
\owncount
\xi\intomega|\nabla c|^2 |u|^2\leq \xi\|\nabla c\|_{\Lo[\frac{2q}{q-2}]}^2\|u \|_{\Lo[q]}^2
\end{align}
is valid for all $t\in(0,\Tm)$. An application of the \GNI\ combined with \cite[Theorem 3.4]{sima90m} allows us to further estimate
\begin{align*}
\|\nabla c\|_{\Lo[\frac{2q}{q-2}]}^2&\leq C_1\|\Delta c\|_{\Lo[2]}^{\frac{4q+4}{3q}}\| c\|_{\Lo[1]}^{\frac{2q-4}{3q}}+C_1\|c\|_{\Lo[1]}^2\\
&\leq C_2\|\Delta c\|_{\Lo[2]}^{\frac{4}{3}+\frac{4}{3q}}+C_2\quad\mbox{for all }t\in(0,\Tm)
\end{align*}
for some $C_1>0$ and $C_2>0$ in view of \eqref{L1-bound-c}. Plugging this into \eqref{c-energy-eq1} and recalling Proposition \ref{Prop:u_bound_N2_all_p}, we thus find $C_3>0$ such that
\begin{align*}
\xi\intomega|\nabla c|^2|u|^2\leq C_3\|\Delta c\|_{\Lo[2]}^{\frac{4}{3}+\frac{4}{3q}}+C_3\quad\mbox{for all }t\in(0,\Tm).
\end{align*}
Since $q>2$, we have $\frac{4}{3}+\frac{4}{3q}<2$ and may apply Young's inequality to obtain 
\begin{align}\label{c-energy-eq2}
\owncount
\xi\intomega|\nabla c|^2|u|^2\leq \frac{\xi}{4}\|\Delta c\|_{\Lo[2]}^2+C_4,
\end{align}
for some $C_4>0$ and all $t\in(0,\Tm)$. To estimate the term containing $f(n)^{2}$ in \eqref{c-energy-eq0} we let $K_1$ denote the positive constant from Corollary \ref{cor:ngnb}. Then, recalling \eqref{fprop} and making use of the fact $\alpha<1$, an application of Young's inequality yields $C_5>0$ fulfilling $\xi f(n)^{2}\leq \frac{1}{2K_1}n^2+C_5$ for all $(x,t)\in\Omega\times(0,\Tm)$ and thus, by Corollary \ref{cor:ngnb}
\begin{align}\label{c-energy-eq3}
\owncount
\xi\intomega f(n)^{2}\leq \frac{1}{2K_1}\intomega n^2+C_5|\Omega|\leq \frac{1}{2}\intomega|\nabla(n^\nfrac{1}{2})|^2+C_6\quad\mbox{for all }t\in(0,\Tm)
\end{align}
with $C_6:=\frac{1}{2}+C_5|\Omega|$. Combining \eqref{c-energy-eq0}, \eqref{c-energy-eq2} and \eqref{c-energy-eq3} completes the proof.
\end{bew}

Before we are able to combine the previous lemmata to derive an ODI appropriate for our purpose, we require one additional result which is a corollary from Lemma \ref{Lem:ngnb}.

\begin{cor}\label{cor:ln_n-gradbound}
There exists a constant $K_4>0$ such that the solution to \eqref{KSaS} obeys
\begin{align*}
\frac{1}{2}\intomega|\nabla(n^{\nfrac{1}{2}})|^2\geq K_4\intomega n\ln n -\frac{1}{2}\quad \mbox{for all }t\in(0,\Tm). 
\end{align*}
\end{cor}

\begin{bew}
In view of the pointwise inequality $x\ln x\leq x^2$ for $x\in(0,\infty)$, the positivity of $n$ ascertained in Lemma \ref{Lem:locEx} therefore implies $n\ln n\leq n^2$ for all $t\in(0,\Tm)$ and thus an application of Corollary \ref{cor:ngnb} immediately shows that there exists $C_1>0$ such that
\begin{align*}
\intomega n\ln n\leq\intomega n^2\leq C_1\|\nabla(n^{\nfrac{1}{2}})\|_{\Lo[2]}^2+C_1
\end{align*}
holds for all $t\in(0,\Tm)$. Therefore, multiplying by $K_4:=\frac{1}{2C_1}$ and reordering the terms appropriately proves the asserted inequality. 
\end{bew}

Adding up suitable multiples of the differential inequalities in Lemma \ref{Lem:n-energy} and Lemma \ref{Lem:c-energy}, we obtain a first bound on the gradient of $c$.

\begin{prop}\label{Prop:grad_c2-bound}
There exists a constant $C>0$ such that the solution of \eqref{KSaS} fulfills
\begin{align}\label{grad_c2-bound}
\owncount
\intomega|\nabla c|^2\leq C\quad\mbox{for all }t\in(0,\Tm).
\end{align}
\end{prop}

\begin{bew}
Letting $K_2$ denote the positive constant from Lemma \ref{Lem:n-energy}, we set $\xi=4K_2+4$ and then $K_3>0$ as the corresponding constant given by Lemma \ref{Lem:c-energy}. With the constants defined this way, we know that the inequality
\begin{align}\label{grad_c2_sptemp-eq1}
\owncount
(2K_2+2)\frac{\intd}{\intdt}\intomega|\nabla c|^2+(K_2+1)\intomega|\Delta c|^2+(4K_2+4)\intomega|\nabla c|^2\leq\frac{1}{2}\intomega|\nabla(n^\nfrac{1}{2})|^2+K_3,
\end{align}
holds for all $t\in(0,\Tm)$ due to Lemma \ref{Lem:c-energy}. Thus, adding up \eqref{n-energy} and \eqref{grad_c2_sptemp-eq1} we obtain
\begin{align*}
\frac{\intd}{\intdt}\bigg(\intomega n\ln n+(2K_2+2)\intomega|\nabla c|^2\bigg)+\frac{1}{2}\intomega|\nabla(n^\nfrac{1}{2})|^2+\intomega|\Delta c|^2+(4K_2+4)\intomega|\nabla c|^2\leq C_1
\end{align*}
for all $t\in(0,\Tm)$ with $C_1=K_2+K_3>0$. By Corollary \ref{cor:ln_n-gradbound} we can estimate $\frac{1}{2}\intomega|\nabla(n^\nfrac{1}{2})|^2$ from below to obtain
\begin{align*}
\frac{\intd}{\intdt}\bigg(\intomega n\ln n+(2K_2+2)\intomega|\nabla c|^2\bigg)+K_4\intomega n\ln n+\!\intomega|\Delta c|^2+2(2K_2+2)\!\intomega|\nabla c|^2\leq C_2
\end{align*}
for all $t\in(0,\Tm)$, with $K_4>0$ as in Corollary \ref{cor:ngnb} and $C_2=C_1+\frac{1}{2}>0$. Dropping the non-negative term involving $|\Delta c|^2$, this implies that $y(t):=\intomega n\ln n+(2K_2+2)\intomega|\nabla c|^2$, $t\in[0,\Tm)$ satisfies
\begin{align*}
y'(t)+C_3y(t)\leq C_2\quad\mbox{for all }t\in(0,\Tm),
\end{align*}
where $C_3:=\min\left\{K_4,2\right\}>0$. Upon an ODE comparison, this leads to the boundedness of $y$ and hence \eqref{grad_c2-bound}, due to $n\ln n$ being bounded from below by the positivity of $n$.
\end{bew}

\subsection{Further testing procedures}\label{ssec:33:testing}

The $L^2$--bound of the gradient of $c$ proven in the previous lemma will be our starting point in improving the regularity of both $n$ and $c$. Preparation and combination of differential inequalities concerning $n^p$ and $|\nabla c|^{2q}$, for appropriately chosen $q$ and $p$, will be the main part of this section. The testing procedures employed in this approach are based on the application to a similar chemotaxis-Stokes system discussed in \cite{win_ct_fluid_3d}. 

The following preparatory result, taken from \cite[Lemma 2.5]{win_ct_fluid_2d}, will be a useful tool in estimations later on and is a simple derivation from Young's inequality.

\begin{lem}\label{young_exp-sum1}
Let $a>0$ and $b>0$ be such that $a+b<1$. Then for all $\varepsilon>0$ there exists $C>0$ such that
\begin{align*}
x^ay^b\leq\varepsilon(x+y)+C\quad\mbox{for all }x\geq0\mbox{ and }y\geq0.
\end{align*}
\end{lem}

In the first step to improve the known regularities of $n$ and $c$, consist of an application of standard testing procedures to gain separate inequalities regarding the time evolution of $\intomega n^p$ and $\intomega|\nabla c|^{2q}$, respectively.

\begin{lem}\label{Lem:np-ineq}
Let $p>1$. Then the solution of \eqref{KSaS} satisfies
\begin{align}\label{np-ineq}
\owncount
\frac{1}{p}\frac{\intd}{\intdt}\intomega n^p+\frac{2(p-1)}{p^2}\intomega|\nabla(n^\nfrac{p}{2})|^2\leq\frac{p-1}{2}\intomega n^p|\nabla c|^2 
\end{align}
for all $t\in(0,\Tm)$.
\end{lem}

\begin{bew}
We multiply the first equation of \eqref{KSaS} with $n^{p-1}$ and integrate by parts to see that
\begin{align*}
\frac{1}{p}\frac{\intd}{\intdt}\intomega n^p=-(p-1)\intomega|\nabla n|^2 n^{p-2}+(p-1)\intomega n^{p-1}\nabla c\cdot\nabla n-\frac{1}{p}\intromega n^{p}u\cdot\vec{\nu}
\end{align*}
holds for all $t\in(0,\Tm)$, where we made use of the fact $\dive u=0$ and the divergence theorem to rewrite the last term accordingly. Due to the boundary condition imposed on $u$ the last term disappears, such that an application of Young's inequality to the second to last term implies
\begin{align*}
\frac{1}{p}\frac{\intd}{\intdt}\intomega n^p+(p-1)\intomega|\nabla n|^2n^{p-2}\leq\frac{p-1}{2}\intomega|\nabla n|^2 n^{p-2}+\frac{p-1}{2}\intomega n^p|\nabla c|^2
\end{align*}
for all $t\in(0,\Tm)$. Reordering the terms and rewriting $|\nabla n|^2n^{p-2}=\frac{4}{p^2}|\nabla(n^\nfrac{p}{2})|^2$ completes the proof.
\end{bew}

\begin{lem}\label{Lem:grad_c2q}
Let $q>1$. Then
\begin{align}\label{grad_c2q}
\owncount
\frac{1}{2q}\frac{\intd}{\intdt}\intomega|\nabla c|^{2q}&+\frac{2(q-1)}{q^2}\intomega\Big\vert\nabla|\nabla c|^q\Big\vert^2+\intomega|\nabla c|^{2q}\nonumber\\
&\leq\left(K_0(q-1)+\frac{K_0}{\sqrt{2}}\right)^2\intomega n^{2\alpha}|\nabla c|^{2q-2}+\intomega|\nabla c|^{2q}|Du|
\end{align}
for all $t\in(0,\Tm)$.
\end{lem}

\begin{bew}
Differentiating the second equation of \eqref{KSaS} and making use of the fact $\Delta|\nabla c|^2=2\nabla c\cdot \nabla\Delta c+2|D^2c|^2$, we obtain
\begin{align*}
\frac{1}{2}\left(|\nabla c|^2\right)_t&=\nabla c\cdot\nabla\left(\Delta c-c+f(n)-u\cdot\nabla c\right)\\
&=\frac{1}{2}\Delta|\nabla c|^2-|D^2c|^2-|\nabla c|^2+\nabla c\cdot\nabla f(n)-\nabla c\cdot\nabla\left(u\cdot\nabla c\right)\mbox{ in }\Omega\times(0,\Tm).
\end{align*} 
We multiply this identity by $\left(|\nabla c|^2\right)^{q-1}$ and integrate by parts over $\Omega$, where due to the Neumann boundary conditions imposed on $n$ and $c$ every boundary integral except the one involving $\frac{\partial|\nabla c|^2}{\partial \nu}$ disappears. Thus, the equality
\begin{align}\label{grad_c2q-eq1}
\owncount
\frac{1}{2q}\frac{\intd}{\intdt}\intomega|\nabla c|^{2q}&+\frac{q-1}{2}\intomega|\nabla c|^{2q-4}\Big\vert\nabla|\nabla c|^2\Big\vert^2+\intomega|\nabla c|^{2q-2}|D^2 c|^2+\intomega|\nabla c|^{2q}\\
&=\intomega|\nabla c|^{2q-2}\nabla c\cdot\nabla f(n)-\intomega|\nabla c|^{2q-2}\nabla c\cdot\nabla\left(u\cdot\nabla c\right)+\frac{1}{2}\intromega|\nabla c|^{2q-2}\frac{\partial|\nabla c|^2}{\partial\nu}\nonumber
\end{align}
holds for all $t\in(0,\Tm)$. Recalling \eqref{fprop}, we integrate the first integral by parts to see that
\begin{align*}
\intomega|\nabla c|^{2q-2}\nabla c\cdot\nabla f\left(n\right)\leq K_0\intomega\Big\vert\nabla|\nabla c|^{2q-2}\Big\vert|\nabla c|\,n^\alpha +K_0\intomega|\nabla c|^{2q-2}|\Delta c|\,n^\alpha
\end{align*}
holds for all $t\in(0,\Tm)$. Since $\nabla|\nabla c|^{2q-2}=2(q-1)|\nabla c|^{2q-4} D^2 c\cdot\nabla c$ in $\Omega\times(0,\Tm)$, and since the Cauchy-Schwarz inequality implies $|\Delta c|\leq \sqrt{2}|D^2 c|$, we may apply Young's inequality to obtain
\begin{align}\label{grad_c2q-eq2}
\owncount
\intomega|\nabla c|^{2q-2}\nabla c\cdot\nabla f\left(n\right)&\leq\intomega|\nabla c|^{2q-2}|D^2 c|^2+\frac{\left(2K_0(q-1)+\sqrt{2}K_0\right)^2}{4}\intomega|\nabla c|^{2q-2} n^{2\alpha}\nonumber\\
&=\intomega|\nabla c|^{2q-2}|D^2 c|^2+\left(K_0(q-1)+\frac{K_0}{\sqrt{2}}\right)^2\intomega|\nabla c|^{2q-2} n^{2\alpha}
\end{align}
for all $t\in(0,\Tm)$. To treat the second integral on the right hand side of \eqref{grad_c2q-eq1}, we first rewrite
\begin{align}\label{grad_c2q-eq3}
\owncount
-\intomega|\nabla c|^{2q-2}\nabla c\cdot\nabla\left(u\cdot \nabla c\right)&=-\intomega|\nabla c|^{2q-2}\nabla c\cdot\left(Du\cdot\nabla c\right)-\intomega|\nabla c|^{2q-2}\nabla c\cdot\left(D^2 c\cdot u\right)
\end{align}
for all $t\in(0,\Tm)$, and then make use of the pointwise equality
\begin{align*}
|\nabla c|^{2q-2}\nabla c\cdot\left(D^2c\cdot u\right)=\frac{1}{2q}u\cdot\nabla|\nabla c|^{2q}\mbox{ in }\Omega\times(0,\Tm),
\end{align*}
to see that, since $u$ is divergence free,
\begin{align*}
-\intomega|\nabla c|^{2q-2}\nabla c\cdot\left(D^2 c\cdot u\right)=\frac{1}{2q}\intomega(\dive u)|\nabla c|^{2q}=0
\end{align*}
holds for all $t\in(0,\Tm)$. Thus, \eqref{grad_c2q-eq3} implies
\begin{align}\label{grad_c2q-eq4}
\owncount
-\intomega|\nabla c|^{2q-2}\nabla c\cdot\nabla\left(u\cdot \nabla c\right)\leq\intomega|\nabla c|^{2q}|Du|\quad\mbox{for all }t\in(0,\Tm).
\end{align}
For the remaining boundary integral in \eqref{grad_c2q-eq1} we recall that the convexity of $\Omega$ ensures $\frac{\partial|\nabla c|^2}{\partial\nu}\leq0$ on $\romega$ (see \cite[Lemme I.1, p.350]{lion}). Combining this with \eqref{grad_c2q-eq1}, \eqref{grad_c2q-eq2} and \eqref{grad_c2q-eq4} completes the proof due to the identity
\[
\Big\vert\nabla|\nabla c|^q\Big\vert^2=\Big\vert\nabla\left(|\nabla c|^{2\cdot\nfrac{q}{2}}\right)\Big\vert^2=\frac{q^2}{4}|\nabla c|^{2q-4}\Big\vert\nabla|\nabla c|^2\Big\vert^2\quad\mbox{in }\Omega\times(0,\Tm).\quad\mbox{\hfill}\qedhere
\]
\end{bew}

Before uniting the inequalities from \eqref{np-ineq} and \eqref{grad_c2q} into a single energy-type inequality, we estimate the right hand sides therein separately.

\begin{lem}\label{Lem:right_hand_sides}
Let $\infty>q>\max\{2,\frac{1}{\alpha}\}$, $p=\alpha q $. For any $\kappa>0$ there exist constants $K_5,K_6$ and $K_7>0$ such that
\begin{align}\label{right_hand_sides-1}
\owncount
\frac{p-1}{2}\intomega n^p|\nabla c|^2\leq \frac{\kappa}{6}\left(\intomega|\nabla(n^\nfrac{p}{2})|^2+\intomega\Big\vert\nabla|\nabla c|^q\Big\vert^2\right)+K_5,
\end{align}
\begin{align}\label{right_hand_sides-2}
\owncount
\left(K_0(q-1)+\frac{K_0}{\sqrt{2}}\right)^2\!\intomega n^{2\alpha}|\nabla c|^{2q-2}\leq\frac{\kappa}{6}\left(\intomega|\nabla(n^\nfrac{p}{2})|^2+\intomega\Big\vert\nabla|\nabla c|^q\Big\vert^2\right)+K_6
\end{align}
and
\begin{align}\label{right_hand_sides-3}
\owncount
\intomega |\nabla c|^{2q}|Du|\leq\frac{\kappa}{6}\intomega\Big\vert\nabla|\nabla c|^q\Big\vert^2+K_7
\end{align}
hold for all $t\in(0,\Tm)$.
\end{lem}

\begin{bew}
To prove \eqref{right_hand_sides-1}, we first fix some $\beta_1>1$ and apply Hölder's inequality to obtain
\begin{align}\label{right_hand_sides-eq1}
\owncount
\frac{p-1}{2}\intomega n^p|\nabla c|^2\leq\frac{p-1}{2}\left(\intomega n^{p\beta_1}\right)^{\nfrac{1}{\beta_1}}\left(\intomega |\nabla c|^{2\beta'_1}\right)^{\nfrac{1}{\beta'_1}}
\end{align}
for all $t\in(0,\Tm)$, where $\beta'_1$ denotes the Hölder conjugate of $\beta_1$. By \eqref{masscons-n} and Lemma \ref{Lem:ngnb} applied to $\varphi=n$, $\eta=m$, $r=p$ and $s=\beta_1$, we can find $C_1>0$ such that
\begin{align}\label{right_hand_sides-eq1.3}
\owncount
\left(\intomega n^{p\beta_1}\right)^\nfrac{1}{\beta_1}\leq C_1\left(\intomega|\nabla(n^\nfrac{p}{2})|^2\right)^{1-\nfrac{1}{p\beta_1}}+C_1\quad\mbox{for all }t\in(0,\Tm).
\end{align}
Similarly to the application of the \GNI\ (\cite[Lemma 2.3]{lankchapto15}) utilized in Lemma \ref{Lem:ngnb}, we can show that the second integral on the right in \eqref{right_hand_sides-eq1} satisfies
\begin{align}\label{right_hand_sides-eq1.6}
\owncount
\left(\intomega |\nabla c|^{2\beta'_1}\right)^{\nfrac{1}{\beta'_1}}\leq C_2\left(\Big\|\nabla|\nabla c|^q\Big\|_{\Lo[2]}^\nfrac{2b_1}{q}\Big\||\nabla c|^q\Big\|_{\Lo[\nfrac{2}{q}]}^\nfrac{(2-2b_1)}{q}+\Big\||\nabla c|^q\Big\|^{\nfrac{2}{q}}_{\Lo[\nfrac{2}{q}]}\right)
\end{align}
for all $t\in(0,\Tm)$ with $C_2>0$ and $b_1\in(0,1)$ provided by
\begin{align*}
b_1=\frac{\frac{q}{2}-\frac{q}{2\beta'_1}}{\frac{q}{2}+\frac{1}{2}-\frac{1}{2}}=1-\frac{1}{\beta'_1}=\frac{1}{\beta_1}.
\end{align*}
Since Proposition \ref{Prop:grad_c2-bound} implies the boundedness of $\||\nabla c|^q\|_{\Lo[\nfrac{2}{q}]}$, plugging \eqref{right_hand_sides-eq1.3} and \eqref{right_hand_sides-eq1.6} into \eqref{right_hand_sides-eq1} we obtain $C_3>0$ such that
\begin{align*}
\frac{p-1}{2}\intomega n^p|\nabla c|^2
\leq C_3&\bigg(\intomega|\nabla(n^\nfrac{p}{2})|^2\bigg)^{1-\nfrac{1}{p\beta_1}}\left(\intomega\Big\vert\nabla|\nabla c|^q\Big|^2\right)^{\nfrac{1}{q\beta_1}}\\
+&\ C_3\left(\intomega|\nabla(n^\nfrac{p}{2})|^2\right)^{1-\nfrac{1}{p\beta_1}}+C_3\left(\intomega\Big\vert\nabla|\nabla c|^q\Big|^2\right)^{\nfrac{1}{q\beta_1}}+C_3
\end{align*}
holds for all $t\in(0,\Tm)$. Due to $\alpha<1$ the choice of $p=\alpha q$					
 implies $p<q$ and thus, $1-\frac{1}{p\beta_1}+\frac{1}{q\beta_1}<1$. Therefore, we may apply Lemma \ref{young_exp-sum1} with $\varepsilon=\frac{\kappa}{12}$ to the three terms on the right hand side containing an integral and obtain for some $C_4>0$ that
\begin{align*}
\frac{p-1}{2}\intomega n^p|\nabla c|^2\leq\frac{\kappa}{6}\left(\intomega|\nabla(n^\nfrac{p}{2})|^2+\intomega\Big\vert\nabla|\nabla c|^q\Big|^2\right)+C_4
\end{align*}
holds for all $t\in(0,\Tm)$, which proves \eqref{right_hand_sides-1}. The proof of \eqref{right_hand_sides-2} follows the same reasoning. First, we apply Hölder's inequality with $\beta_2=\frac{q+1}{2}$ and $\beta'_2$ as corresponding Hölder conjugate
to obtain
\begin{align}\label{right_hand_sides-eq2}
\owncount
\intomega n^{2\alpha}|\nabla c|^{2q-2}\leq\left(\intomega n^{2\alpha\beta_2}\right)^{\nfrac{1}{\beta_2}}\left(\intomega |\nabla c|^{(2q-2)\beta'_2}\right)^{\nfrac{1}{\beta'_2}}
\end{align}
for all $t\in(0,\Tm)$. Since the choices of $\beta_2$ and $p$ imply $\frac{2\alpha\beta_2}{p}=\frac{\alpha(q+1)}{\alpha q}>1$, we can utilize Lemma \ref{Lem:ngnb} with $\varphi=n$, $r=p$ and $s=\frac{2\alpha\beta_2}{p}$ to estimate
\begin{align}\label{right_hand_sides-eq3}
\owncount
\left(\intomega n^{2\alpha\beta_2}\right)^{\nfrac{1}{\beta_2}}\leq C_5\left(\intomega|\nabla(n^\nfrac{p}{2})|^2\right)^{\nfrac{(2\alpha\beta_2-1)}{p\beta_2}}+C_5\quad\mbox{for all }t\in(0,\Tm),
\end{align}
with some $C_5>0$. For the integral involving $|\nabla c|^{(2q-2)\beta'_2}$, we make use of the \GNI\ as shown before to obtain $C_6>0$ such that
\begin{align}\label{right_hand_sides-eq4}
\owncount
\left(\intomega |\nabla c|^{(2q-2)\beta'_2}\right)^{\nfrac{1}{\beta'_2}}\leq C_6\left(\intomega\Big|\nabla|\nabla c|^q\Big|^2\right)^{\frac{q-1}{q}b_2}+C_6
\end{align}
holds for all $t\in(0,\Tm)$, with $b_2\in(0,1)$ determined by
\begin{align*}
b_2=\frac{\frac{q}{2}-\frac{q}{2(q-1)\beta'_2}}{\frac{q}{2}+\frac{1}{2}-\frac{1}{2}}=1-\frac{1}{(q-1)\beta'_2}=1-\frac{1}{(q-1)}+\frac{1}{(q-1)\beta_2}.
\end{align*}
Thus, a combination of \eqref{right_hand_sides-eq2},\eqref{right_hand_sides-eq3} and \eqref{right_hand_sides-eq4} leads to
\begin{align*}
\bigg(K_0(q-1)+\frac{K_0}{\sqrt{2}}&\bigg)^2\!\intomega n^{2\alpha}|\nabla c|^{2q-2}\leq C_7\bigg(\intomega|\nabla(n^\nfrac{p}{2})|^2\bigg)^{\nfrac{(2\alpha\beta_2-1)}{p\beta_2}}\left(\intomega\Big|\nabla|\nabla c|^q\Big|^2\right)^{\frac{q-1}{q}b_2}\\
&\ +C_7\left(\intomega|\nabla(n^\nfrac{p}{2})|^2\right)^{\nfrac{(2\alpha\beta_2-1)}{p\beta_2}}+C_7\left(\intomega\Big|\nabla|\nabla c|^q\Big|^2\right)^{\frac{q-1}{q}b_2}+C_7
\end{align*}
for all $t\in(0,\Tm)$ with some $C_7>0$.
Here the choice of $p$ and the fact that $\alpha<1$ imply
\begin{align*}
\frac{2\alpha\beta_2-1}{p\beta_2}+\frac{q-1}{q}b_2&=\frac{2\alpha}{p}-\frac{1}{p\beta_2}+\frac{q-2}{q}+\frac{1}{q\beta_2}\\
&=\frac{2}{q}-\frac{1}{\alpha q\beta_2 }+\frac{q-2}{q}+\frac{1}{q\beta_2 }=1-\frac{1-\alpha}{\alpha q\beta_2 }<1.
\end{align*}
Therefore, the requirements of Lemma \ref{young_exp-sum1} are satisfied again and an application thereof yields $C_8>0$ such that
\begin{align*}
\left(K_0(q-1)+\frac{K_0}{\sqrt{2}}\right)^2\!\intomega n^{2\alpha}|\nabla c|^{2q-2}\leq \frac{\kappa}{6}\left(\intomega|\nabla(n^\nfrac{p}{2})|^2+\intomega\Big|\nabla|\nabla c|^q\Big|^2\right)+C_8
\end{align*}
holds for all $t\in(0,\Tm)$ and thus proving \eqref{right_hand_sides-2}. To verify \eqref{right_hand_sides-3} we fix $\beta_3=\frac{3}{2}$ and $\beta'_3=3$. Since $\beta_3<2$ Hölder's inequality yields
\begin{align*}
\intomega|\nabla c|^{2q}|Du|\leq\left(\intomega|\nabla c|^{2q\beta'_3}\right)^{\nfrac{1}{\beta'_3}}\left(\intomega|D u|^{\beta_3}\right)^{\nfrac{1}{\beta_3}}\leq C_9\Big\||\nabla c|^q\Big\|_{\Lo[6]}^2
\end{align*}
for some $C_9>0$, in view of the boundedness of $\|D u\|_{\Lo[\frac{3}{2}]}$ shown in Proposition \ref{Prop:u_bound_N2_all_p}. Similarly to the previous applications of the Gagliardo-Nirenberg and Young inequalities we can make use of the boundedness of $\||\nabla c|^q\|_{\Lo[\nfrac{2}{q}]}$ to obtain $C_{10}>0$ such that
\begin{align*}
\intomega|\nabla c|^{2q}|Du|\leq\frac{\kappa}{6}\intomega\Big|\nabla|\nabla c|^q\Big|^2+C_{10}
\end{align*}
for all $t\in(0,\Tm)$, which completes the proof.
\end{bew}

Combining the three previous lemmata we are now in the position to control $L^p$--norms of $n$ and  $\nabla c$ with arbitrarily high $p$. In fact we have

\begin{prop}\label{Prop:np_grad_c2q-bounds}
Let $\infty>q>\max\{2,\frac{1}{\alpha}\}$ and $p=\alpha q$. Then we can find $C>0$ such that, the solution to \eqref{KSaS} satisfies
\begin{align}\label{np-bound}
\owncount
\intomega n^p\leq C\quad\mbox{for all }t\in(0,\Tm)
\end{align}
and
\begin{align}\label{grad_c2q-bound}
\owncount
\intomega |\nabla c|^{2q}\leq C\quad\mbox{for all }t\in(0,\Tm).
\end{align}
\end{prop}

\begin{bew}
Given $q>\max\{2,\frac{1}{\alpha}\}$ and $p=\alpha q$ we fix $\kappa=\min\left\{\frac{2(q-1)}{q^2},\,\frac{2(p-1)}{p^2}\right\}$. By Lemmata \ref{Lem:np-ineq}, \ref{Lem:grad_c2q} and \ref{Lem:right_hand_sides}, we can find $C_1:=K_5+K_6+K_7>0$ such that
\begin{align*}
\frac{\intd}{\intdt}\bigg(\frac{1}{p}\intomega n^p&+\frac{1}{2q}\intomega|\nabla c|^{2q}\bigg)+\frac{2(p-1)}{p^2}\intomega|\nabla(n^\nfrac{p}{2})|^2\\
&+\frac{2(q-1)}{q^2}\intomega\Big\vert\nabla|\nabla c|^q\Big\vert^2+\intomega|\nabla c|^{2q}
\leq\frac{\kappa}{2}\left(\intomega|\nabla(n^\nfrac{p}{2})|^2+\intomega\Big\vert\nabla|\nabla c|^q\Big\vert^2\right)+C_1
\end{align*}
holds for all $t\in(0,\Tm)$. Herein the choice of $\kappa$ implies
\begin{align}\label{np_grad_c2q-bounds-eq1}
\owncount
\frac{\intd}{\intdt}\bigg(\frac{1}{p}\intomega n^p+\frac{1}{2q}\intomega|\nabla c|^{2q}\bigg)+\frac{p-1}{p^2}\intomega|\nabla(n^\nfrac{p}{2})|^2+\frac{q-1}{q^2}\intomega\Big\vert\nabla|\nabla c|^q\Big\vert^2+\intomega|\nabla c|^{2q}\leq C_1
\end{align}
for all $t\in(0,\Tm)$. We drop the non-negative term $\frac{q-1}{q^2}\intomega|\nabla|\nabla c|^q|^2$ and apply Lemma \ref{Lem:ngnb} to estimate $\intomega|\nabla(n^\nfrac{p}{2})|^2$ from below in \eqref{np_grad_c2q-bounds-eq1}, to obtain $C_2,C_3>0$ such that $y(t):=\frac{1}{p}\intomega n^p+\frac{1}{2q}\intomega|\nabla c|^{2q}$, $t\in(0,\Tm)$ satisfies
\begin{align*}
y'(t)+C_2 y(t)\leq C_3\quad\mbox{for all }t\in(0,\Tm),
\end{align*}
from which we infer the boundedness of $y$ upon an ODE comparison and thus \eqref{np-bound} and \eqref{grad_c2q-bound}.
\end{bew}

\subsection{Global existence and boundedness}\label{ssec:34:bound}
We can now begin to verify the boundedness of the three quantities appearing in the extensibility criterion \eqref{lExAlt}. The first of these quantities will be $\|A^\delta u(\cdot,t)\|_{\Lo[2]}$.

\begin{prop}\label{Prop:Au2-bound}
Let $\delta\in(\frac{1}{2},1)$ be as in Lemma \ref{Lem:locEx}. There exists a constant $C>0$ such that the solution of \eqref{KSaS} satisfies
\begin{align*}
\|A^\delta u(\cdot,t)\|_{\Lo[2]}\leq C\quad\mbox{for all }t\in(0,\Tm).
\end{align*}
\end{prop}

\begin{bew}
 The proof essentially follows the argumentation of \cite[Lemma 2.3]{win_ct_fluid_2d}, whilst making use of the previously proven bound $\|n\|_{\Lo[p]}\leq C$ for all $t\in(0,\Tm)$ with some $p>2$. Nonetheless, let us recount the main arguments.

It is well known, see \cite[Theorem 38.6]{sellyou} and \cite[p.204]{sohr} for instance, that the Stokes operator $A$ is a positive, sectorial operator and generates a contraction semigroup $\left(\etA\right)_{t\geq0}$ in $L^2_\sigma(\Omega)$ with operator norm bounded by
\begin{align*}
\|\etA\|\leq e^{-\mu_1 t}\quad\mbox{for all }t\geq0,
\end{align*}
with some $\mu_1>0$. Furthermore, the operator norm of the fractional powers of the Stokes operator satisfy an exponential decay property (\cite[Theorem 37.5]{sellyou}). That is, there exists $C_1>0$ such that
\begin{align}\label{Au2-bound-expdecprop}
\owncount
\left\| A^\delta\etA\right\|\leq C_1t^{-\delta}e^{-\mu_1 t}\quad\mbox{for all }t>0.
\end{align}
Thus, representing $u$ by its variation of constants formula
\begin{align*}
u(\cdot,t)=\etA u_0+\intot e^{-(t-s)A}\Po\left(n(\cdot,s)\nabla\phi\right)\intd s,\quad t\in(0,\Tm),
\end{align*}
and applying the fractional power $A^\delta$, we can make use of the fact that $\etA$ commutes with $A^\delta$ (\cite[IV.(1.5.16), p.206]{sohr}), the contraction property and \eqref{Au2-bound-expdecprop} to find $C_2>0$ such that
\begin{align}\label{Au2-bound-eq1}
\owncount
\|A^\delta u(\cdot,t)\|_{\Lo[2]}&\leq \|A^\delta u_0\|_{\Lo[2]}+C_1\intot \left(t-s\right)^{-\delta}e^{-\mu_1(t-s)}\left\|\Po\left(n(\cdot,s)\nabla\phi\right)\right\|_{\Lo[2]}\intd s\nonumber\\
&\leq \|A^\delta u_0\|_{\Lo[2]}+C_2\sup_{t\in(0,\Tm)}\left\|n(\cdot,t)\right\|_{\Lo[2]}\int_0^\infty\!\! \sigma^{-\delta} e^{-\mu_1\sigma}\intd \sigma
\end{align}
holds for all $t\in(0,\Tm)$, by the boundedness of $\nabla\phi$. Due to \eqref{idreg} we have $\|A^\delta u_0\|_{\Lo[2]}\leq C_3$ for some $C_3>0$. Furthermore, since $\delta<1$ the integral converges and by Proposition \ref{Prop:np_grad_c2q-bounds} we can find $C_4>0$ such that $\|n(\cdot,t)\|_{\Lo[2]}\leq C_4$ for all $t\in(0,\Tm)$. Combined with \eqref{Au2-bound-eq1} these facts yield
\begin{align*}
\|A^\delta u(\cdot,t)\|_{\Lo[2]}\leq C_5\quad\mbox{for all }t\in(0,\Tm)
\end{align*}
with some $C_5>0$, which completes the proof.
\end{bew}

The second quantity of the extensibility criterion we treat is $\|c(\cdot,t)\|_{\W[1,\theta]}$. In view of Proposition  \ref{Prop:np_grad_c2q-bounds}, we can take some $q>\max\{2,\theta\}$ and obtain under simple application of the Poincaré inequality 

\begin{cor}\label{Cor:grad_ctheta-bound}
There exists a constant $C>0$ such that
\begin{align*}
\|c(\cdot,t)\|_{\W[1,\theta]}\leq C
\end{align*}
holds for all $t\in(0,\Tm)$.
\end{cor}

Now, to prove the last remaining bound required for the extensibility criterion \eqref{lExAlt}, as well as one of the estimates required for the boundedness result, we require some well known results concerning the Neumann heat semigroup $\left(\etd\right)_{t\geq0}$. These semigroup estimates and Proposition \ref{Prop:np_grad_c2q-bounds} will be the main ingredients of our proof. For more details concerning the estimations used, we refer the reader to \cite[Lemma 2.1]{cao2014global}, \cite[Lemma 1.3]{win10jde} and \cite{hen81}.

\begin{prop}\label{Prop:ninf-bound}
There exists a constant $C>0$ such that
\begin{align*}
\|n(\cdot,t)\|_{\Lo[\infty]}\leq C
\end{align*}
holds for all $t\in(0,\Tm)$.
\end{prop}

\begin{bew}
First, we fix $p>2$ and represent $n$ by its variation of constants formula
\begin{align*}
n(\cdot,t)=\etd n_0-\intot\etsd\big(\nabla\cdot\left(n\nabla c\right)+u\cdot\nabla n\big)(\cdot,s)\intd s,\quad t\in(0,\Tm).
\end{align*}
The fact $\dive u=0$ and the maximum principle therefore yield
\begin{align*}
\|n(\cdot,t)\|_{\Lo[\infty]}&\leq \| n_0\|_{\Lo[\infty]}+\intot\left\|\etsd\big(\nabla\cdot\left(n\nabla c+u n\right)\big)(\cdot,s)\right\|_{\Lo[\infty]}\intd s
\end{align*}
for all $t\in(0,\Tm)$. Now, we can make use of the well known smoothing properties of the Neumann heat semigroup (see \cite[Lemma 2.1  (iv)]{cao2014global}), to estimate
\begin{align}\label{ninf-bound-eq1}
\owncount
\|n(\cdot,t)\|_{\Lo[\infty]}&\leq \|n_0\|_{\Lo[\infty]}\\&\ +C_1\int_0^{t}\!\!\left(1+\left(t-s\right)^{-\frac{1}{2}-\frac{1}{p}}\right)e^{-\lambda_1(t-s)}\left(\left\|\left(n\left(\nabla c +u\right)\right)(\cdot,s)\right\|_{\Lo[p]}\right)\intd s\nonumber
\end{align}
for all $t\in(0,\Tm)$ and some $C_1>0$, where $\lambda_1$ denotes the first nonzero eigenvalue of $-\Delta$ in $\Omega$ with regards to the homogeneous Neumann boundary conditions. To estimate $\|n(\nabla c+u)\|_{\Lo[p]}$ we apply Hölder's inequality to obtain some $C_2>0$ such that
\begin{align*}
\|n(\nabla c+u)(\cdot,t)\|_{\Lo[p]}\leq\|n(\cdot,t)\|_{\Lo[2p]}\left(\|\nabla c(\cdot,t)\|_{\Lo[2p]}+\|u(\cdot,t)\|_{\Lo[2p]}\right)\leq C_2
\end{align*}
holds for all $t\in(0,\Tm)$, wherein the boundedness of all quantities on the right hand side followed in view of Propositions \ref{Prop:u_bound_N2_all_p} and \ref{Prop:np_grad_c2q-bounds}. Plugging this into \eqref{ninf-bound-eq1} and recalling $n_0\in \CSp{0}{\bomega}$ yields $C_3>0$ such that
\begin{align*}
\|n(\cdot,t)\|_{\Lo[\infty]}&\leq C_3+C_3\int_0^\infty\!\!\left(1+\left(t-s\right)^{-\frac{1}{2}-\frac{1}{p}}\right)e^{-\lambda_1(t-s)}\intd s
\end{align*}
is valid for all $t\in(0,\Tm)$. By the choice of $p$ we have $-\frac{1}{2}-\frac{1}{p}>-1$ and thus there exists $C_4>0$ such that
\begin{align*}
\|n(\cdot,t)\|_{\Lo[\infty]}\leq C_4\quad\mbox{for all }t\in(0,\Tm),
\end{align*}
which completes the proof.
\end{bew}

Let us gather the previous results to prove our main theorem.

\begin{proof}[\textbf{Proof of Theorem \ref{Thm:globEx}:}]
As an immediate consequence of the bounds in Proposition \ref{Prop:Au2-bound}, Corollary \ref{Cor:grad_ctheta-bound} and Proposition \ref{Prop:ninf-bound}, we obtain $\Tm=\infty$ in view of the extensibility criterion \eqref{lExAlt}.
Secondly, since $\theta>2$ we have $\W[1,\theta]\hookrightarrow\Co[\mu_1]$ with $\mu_1=\frac{\theta-2}{\theta}$ (\cite[Theorem 5.6.5]{evans}). Thus, Corollary \ref{Cor:grad_ctheta-bound} implies $\|c(\cdot,t)\|_{\Lo[\infty]}\leq C$ for all $t\in(0,\Tm)$. Additionally, since for $\delta\in(\frac{1}{2},1)$ the fractional powers of the Stokes operator satisfy $D(A^\delta)\hookrightarrow\Co[\mu_2]$ for any $\mu_2\in(0,2\delta-1)$ (see \cite[Lemma III.2.4.3]{sohr} and \cite[Theorem 5.6.5]{evans}), Proposition \ref{Prop:Au2-bound} shows that $\|u(\cdot,t)\|_{\Lo[\infty]}\leq C$ for all $t\in(0,\Tm)$ and the boundedness of $\|n(\cdot,t)\|_{\Lo[\infty]}$ for all $t\in(0,\Tm)$ follows directly from Proposition \ref{Prop:ninf-bound}.
\end{proof}



\begin{thebibliography}{45}
\providecommand{\natexlab}[1]{#1}
\providecommand{\url}[1]{\texttt{#1}}
\expandafter\ifx\csname urlstyle\endcsname\relax
  \providecommand{\doi}[1]{doi: #1}\else
  \providecommand{\doi}{doi: \begingroup \urlstyle{rm}\Url}\fi

\bibitem[Bellomo et~al.(2015)Bellomo, Bellouquid, Tao, and Winkler]{BBWT15}
N.~Bellomo, A.~Bellouquid, Y.~Tao, and M.~Winkler.
\newblock Toward a mathematical theory of {K}eller–{S}egel models of pattern
  formation in biological tissues.
\newblock \emph{Math. Mod. Meth. Appl. S.}, 25\penalty0 (09):\penalty0
  1663--1763, 2015.

\bibitem[Cao(2015)]{cao2014global}
X.~Cao.
\newblock Global bounded solutions of the higher-dimensional keller-segel
  system under smallness conditions in optimal spaces.
\newblock \emph{Discrete and Continuous Dynamical Systems}, 35\penalty0
  (5):\penalty0 1891--1904, 2015.

\bibitem[Cao and Lankeit(2016)]{caolan16_smalldatasol3dnavstokes}
X.~Cao and J.~Lankeit.
\newblock Global classical small-data solutions for a three-dimensional
  chemotaxis {N}avier-{S}tokes system involving matrix-valued sensitivities.
\newblock 2016.
\newblock arXiv:1601.03897 - Preprint.

\bibitem[Coll et~al.(1994)Coll, Bowden, Meehan, Konig, Carroll, Tapiolas,
  Alino, Heaton, De~Nys, Leone, et~al.]{coll1994chemical}
J.~Coll, B.~Bowden, G.~Meehan, G.~Konig, A.~Carroll, D.~Tapiolas, P.~Alino,
  A.~Heaton, R.~De~Nys, P.~Leone, et~al.
\newblock Chemical aspects of mass spawning in corals. i. sperm-attractant
  molecules in the eggs of the scleractinian coral montipora digitata.
\newblock \emph{Mar. Biol.}, 118\penalty0 (2):\penalty0 177--182, 1994.

\bibitem[Di~Francesco et~al.(2010)Di~Francesco, Lorz, and
  Markowich]{francescolorz10}
M.~Di~Francesco, A.~Lorz, and P.~Markowich.
\newblock Chemotaxis-fluid coupled model for swimming bacteria with nonlinear
  diffusion: global existence and asymptotic behavior.
\newblock \emph{Discrete Contin. Dyn. Syst.}, 28\penalty0 (4):\penalty0
  1437--1453, 2010.

\bibitem[Evans(2010)]{evans}
L.~C. Evans.
\newblock \emph{Partial differential equations}, volume~19 of \emph{Graduate
  Studies in Mathematics}.
\newblock American Mathematical Society, Providence, RI, second edition, 2010.

\bibitem[Friedman(1969)]{fr69}
A.~Friedman.
\newblock \emph{Partial differential equations}.
\newblock Holt, Rinehart and Winston, 1969.

\bibitem[Giga(1981)]{gig81}
Y.~Giga.
\newblock Analyticity of the semigroup generated by the {S}tokes operator in
  {$L_{r}$} spaces.
\newblock \emph{Math. Z.}, 178\penalty0 (3):\penalty0 297--329, 1981.

\bibitem[Henry(1981)]{hen81}
D.~Henry.
\newblock \emph{Geometric Theory of Semilinear Parabolic Equations}, volume 840
  of \emph{Lecture Notes in Mathematics}.
\newblock Springer Berlin Heidelberg, 1981.

\bibitem[Hillen and Painter(2009)]{HP09}
T.~Hillen and K.~J. Painter.
\newblock A user’s guide to {PDE} models for chemotaxis.
\newblock \emph{J. Math. Biol.}, 58\penalty0 (1-2):\penalty0 183--217, 2009.

\bibitem[Horstmann(2003)]{Ho03}
D.~Horstmann.
\newblock From 1970 until present: the keller-segel model in chemotaxis and its
  consequences i.
\newblock \emph{Jahresber. Deutsch. Math.-Verein.}, 105:\penalty0 103--165,
  2003.

\bibitem[Horstmann and Winkler(2005)]{HoWin05_bvblowchemo}
D.~Horstmann and M.~Winkler.
\newblock Boundedness vs. blow-up in a chemotaxis system.
\newblock \emph{J. Differential Equations}, 215\penalty0 (1):\penalty0 52 --
  107, 2005.

\bibitem[Keller and Segel(1970)]{KS70}
E.~F. Keller and L.~A. Segel.
\newblock Initiation of slime mold aggregation viewed as an instability.
\newblock \emph{J. Theor. Biol.}, 26\penalty0 (3):\penalty0 399--415, 1970.

\bibitem[Keller and Segel(1971)]{KS71}
E.~F. Keller and L.~A. Segel.
\newblock Model for chemotaxis.
\newblock \emph{J. Theor. Biol.}, 30\penalty0 (2):\penalty0 225--234, 1971.

\bibitem[Kozono et~al.(2015)Kozono, Miura, and Sugiyama]{kozono15}
H.~Kozono, M.~Miura, and Y.~Sugiyama.
\newblock Existence and uniqueness theorem on mild solutions to the
  {K}eller-{S}egel system coupled with the {N}avier-{S}tokes fluid.
\newblock 2015.
\newblock Preprint.

\bibitem[Lankeit and Li(2015)]{lankchapto15}
J.~Lankeit and Y.~Li.
\newblock Boundedness in a chemotaxis-haptotaxis model with nonlinear
  diffusion.
\newblock 2015.
\newblock arXiv:1508.05846 - Preprint.

\bibitem[Lions(1980)]{lion}
P.-L. Lions.
\newblock R\'esolution de probl\`emes elliptiques quasilin\'eaires.
\newblock \emph{Arch. Rational Mech. Anal.}, 74\penalty0 (4):\penalty0
  335--353, 1980.

\bibitem[Liu()]{liudongmei15_boundchemo}
D.~Liu.
\newblock Boundedness in a chemotaxis system with nonlinear signal production.
\newblock Preprint.

\bibitem[Lorz(2010)]{lorz10}
A.~Lorz.
\newblock Coupled chemotaxis fluid model.
\newblock \emph{Math. Mod. Meth. Appl. S.}, 20\penalty0 (06):\penalty0
  987--1004, 2010.

\bibitem[Miller(1985)]{miller1985demonstration}
R.~L. Miller.
\newblock Demonstration of sperm chemotaxis in echinodermata: {A}steroidea,
  {H}olothuroidea, {O}phiuroidea.
\newblock \emph{J. Exp. Zool.}, 234\penalty0 (3):\penalty0 383--414, 1985.

\bibitem[Mimura and Tsujikawa(1996)]{Mimura1996499}
M.~Mimura and T.~Tsujikawa.
\newblock Aggregating pattern dynamics in a chemotaxis model including growth.
\newblock \emph{Physica A}, 230\penalty0 (3–4):\penalty0 499 -- 543, 1996.

\bibitem[Mizoguchi and Winkler(2013)]{mizoguchi_winkler_13}
N.~Mizoguchi and M.~Winkler.
\newblock Blow-up in the two-dimensional parabolic {K}eller-{S}egel system.
\newblock 2013.
\newblock Preprint.

\bibitem[Nagai et~al.(1997)Nagai, Senba, and Yoshida]{NSY97}
T.~Nagai, T.~Senba, and K.~Yoshida.
\newblock Application of the {T}rudinger-{M}oser inequality to a parabolic
  system of chemotaxis.
\newblock \emph{Funkcial. Ekvac.}, 40\penalty0 (3):\penalty0 411--433, 1997.

\bibitem[Nakaguchi and Osaki(2013)]{NaOs13}
E.~Nakaguchi and K.~Osaki.
\newblock Global solutions and exponential attractors of a parabolic-parabolic
  system for chemotaxis with subquadratic degradation.
\newblock \emph{Discrete Contin. Dyn. Syst. Ser. B}, 18\penalty0 (10):\penalty0
  2627--2646, 2013.

\bibitem[Osaki and Yagi(2001)]{OY01}
K.~Osaki and A.~Yagi.
\newblock Finite dimensional attractor for one-dimensional keller-segel
  equations.
\newblock \emph{Funkcial. Ekvac.}, 44\penalty0 (3):\penalty0 441--470, 2001.

\bibitem[Osaki et~al.(2002)Osaki, Tsujikawa, Yagi, and
  Mimura]{Os02-chemologatract}
K.~Osaki, T.~Tsujikawa, A.~Yagi, and M.~Mimura.
\newblock Exponential attractor for a chemotaxis-growth system of equations.
\newblock \emph{Nonlinear Anal.}, 51\penalty0 (1, Ser. A: Theory
  Methods):\penalty0 119--144, 2002.

\bibitem[Sell and You(2002)]{sellyou}
G.~R. Sell and Y.~You.
\newblock \emph{Dynamics of evolutionary equations}, volume 143 of
  \emph{Applied Mathematical Sciences}.
\newblock Springer-Verlag, New York, 2002.

\bibitem[Simader(1990)]{sima90m}
C.~G. Simader.
\newblock The weak {D}irichlet and {N}eumann problem for the {L}aplacian in
  {$L^q$} for bounded and exterior domains. {A}pplications.
\newblock In M.~Krbec, A.~Kufner, B.~Opic, and J.~Rákosník, editors,
  \emph{Nonlinear {A}nalysis, {F}unction {S}paces and {A}pplications {V}ol. 4},
  Teubner-Texte Math., pages 180--223. Vieweg+Teubner Verlag, 1990.

\bibitem[Sohr(2001)]{sohr}
H.~Sohr.
\newblock \emph{The {N}avier-{S}tokes equations}.
\newblock Birkh\"auser Advanced Texts: Basler Lehrb\"ucher. Birkh\"auser
  Verlag, Basel, 2001.
\newblock An elementary functional analytic approach.

\bibitem[Tao and Winkler(2011)]{tao_winkler_chemohapto11siam}
Y.~Tao and M.~Winkler.
\newblock A chemotaxis-haptotaxis model: the roles of nonlinear diffusion and
  logistic source.
\newblock \emph{SIAM J. Math. Anal.}, 43\penalty0 (2):\penalty0 685--704, 2011.

\bibitem[Tao and Winkler(2012{\natexlab{a}})]{TaoWin12_evsmooth}
Y.~Tao and M.~Winkler.
\newblock Eventual smoothness and stabilization of large-data solutions in a
  three-dimensional chemotaxis system with consumption of chemoattractant.
\newblock \emph{J. Differential Equations}, 252\penalty0 (3):\penalty0 2520 --
  2543, 2012{\natexlab{a}}.

\bibitem[Tao and Winkler(2012{\natexlab{b}})]{win_ct_fluid_2d}
Y.~Tao and M.~Winkler.
\newblock Global existence and boundedness in a {K}eller-{S}egel-{S}tokes model
  with arbitrary porous medium diffusion.
\newblock \emph{Discrete Contin. Dyn. Syst.}, 32\penalty0 (5):\penalty0
  1901--1914, 2012{\natexlab{b}}.

\bibitem[Tao and Winkler(2015{\natexlab{a}})]{tao_winkler15_zampfinal}
Y.~Tao and M.~Winkler.
\newblock Boundedness and decay enforced by quadratic degradation in a
  three-dimensional chemotaxis--fluid system.
\newblock \emph{Z. Angew. Math. Phys.}, 66\penalty0 (5):\penalty0 2555--2573,
  2015{\natexlab{a}}.

\bibitem[Tao and Winkler(2015{\natexlab{b}})]{tao_winkler_non2015}
Y.~Tao and M.~Winkler.
\newblock Blow-up prevention by quadratic degradation in a two-dimensional
  {K}eller-{S}egel-{N}avier-{S}tokes system.
\newblock 2015{\natexlab{b}}.
\newblock Preprint.

\bibitem[Tuval et~al.(2005)Tuval, Cisneros, Dombrowski, Wolgemuth, Kessler, and
  Goldstein]{tuval2005bacterial}
I.~Tuval, L.~Cisneros, C.~Dombrowski, C.~W. Wolgemuth, J.~O. Kessler, and R.~E.
  Goldstein.
\newblock Bacterial swimming and oxygen transport near contact lines.
\newblock \emph{Proc. Natl. Acad. Sci. U.S.A.}, 102\penalty0 (7):\penalty0
  2277--2282, 2005.

\bibitem[Wang and Xiang(2015)]{Wang20157578}
Y.~Wang and Z.~Xiang.
\newblock Global existence and boundedness in a keller–segel–stokes system
  involving a tensor-valued sensitivity with saturation.
\newblock \emph{J. Differential Equations}, 259\penalty0 (12):\penalty0 7578 --
  7609, 2015.

\bibitem[Wiegner(1999)]{Wie99-NS}
M.~Wiegner.
\newblock The {N}avier-{S}tokes equations---a neverending challenge?
\newblock \emph{Jahresber. Deutsch. Math.-Verein.}, 101\penalty0 (1):\penalty0
  1--25, 1999.

\bibitem[Winkler(2010{\natexlab{a}})]{Win10-chemolog}
M.~Winkler.
\newblock Boundedness in the higher-dimensional parabolic-parabolic chemotaxis
  system with logistic source.
\newblock \emph{Comm. Partial Differential Equations}, 35\penalty0
  (8):\penalty0 1516--1537, 2010{\natexlab{a}}.

\bibitem[Winkler(2010{\natexlab{b}})]{win10jde}
M.~Winkler.
\newblock Aggregation vs. global diffusive behavior in the higher-dimensional
  {K}eller-{S}egel model.
\newblock \emph{J. Differential Equations}, 248\penalty0 (12):\penalty0
  2889--2905, 2010{\natexlab{b}}.

\bibitem[Winkler(2012)]{win_fluid_final}
M.~Winkler.
\newblock Global large-data solutions in a chemotaxis-({N}avier-){S}tokes
  system modeling cellular swimming in fluid drops.
\newblock \emph{Comm. Partial Differential Equations}, 37\penalty0
  (2):\penalty0 319--351, 2012.

\bibitem[Winkler(2015{\natexlab{a}})]{win15_chemorot}
M.~Winkler.
\newblock Large-data global generalized solutions in a chemotaxis system with
  tensor-valued sensitivities.
\newblock \emph{SIAM J. Math. Anal.}, 47\penalty0 (4):\penalty0 3092--3115,
  2015{\natexlab{a}}.

\bibitem[Winkler(2015{\natexlab{b}})]{win15_globweak3d}
M.~Winkler.
\newblock Global weak solutions in a three-dimensional
  chemotaxis-{N}avier-{S}tokes system.
\newblock \emph{Ann. Inst. H. Poincar{é} Anal. Non Lin{é}aire},
  2015{\natexlab{b}}.
\newblock to appear.

\bibitem[Winkler(2015{\natexlab{c}})]{win_ct_fluid_3d}
M.~Winkler.
\newblock Boundedness and large time behavior in a three-dimensional
  chemotaxis-{S}tokes system with nonlinear diffusion and general sensitivity.
\newblock \emph{Calc. Var. Partial Differential Equations}, 54\penalty0
  (4):\penalty0 3789--3828, 2015{\natexlab{c}}.

\bibitem[Winkler(2016)]{win15_chemonavstokesfinal}
M.~Winkler.
\newblock How far do chemotaxis-driven forces influence regularity in the
  {N}avier-{S}tokes system?
\newblock \emph{Trans. Amer. Math. Soc.}, 2016.
\newblock to appear.

\bibitem[Xue and Othmer(2009)]{XO09-MSmodels}
C.~Xue and H.~G. Othmer.
\newblock Multiscale models of taxis-driven patterning in bacterial
  populations.
\newblock \emph{SIAM J. Appl. Math.}, 70\penalty0 (1):\penalty0 133--167, 2009.

\end{thebibliography}

\end{document}